\documentclass[12pt,twoside,a4paper]{amsart}%article
\UseRawInputEncoding
\usepackage{amsmath,amssymb}
\usepackage{cite}
\usepackage{graphicx}
\usepackage{txfonts}
\usepackage{mathrsfs}
\newcommand{\K}{J.Kaczorowski and K.Wiertelak}
\newcommand{\ka}{J.Kaczorowski}

\allowdisplaybreaks[4]
\theoremstyle{theorem}
\newtheorem{theorem}{Theorem}[section]
\newtheorem{lemma}[theorem]{Lemma}
\newtheorem{remark}[theorem]{Remark}

\newtheorem{fact}[theorem]{Fact}

\begin{document}

\begin{center}
\textbf{\large On the Volterra integral equation for the remainder term in the asymptotic formula on the associated Euler totient function}
\\
\vspace{0.9cm}
by
\\
\vspace{0.35cm}
Hideto  IWATA
\end{center}

\vspace{0.5cm}

\textbf{Abstract.}\hspace{0.1cm}{\K} considered the integral equation for remainder terms in the asymptotic formula for the Euler totient function and for the twisted Euler $\varphi$-function. In 2013, {\ka} defined  the associated Euler totient function which extends the above two functions and proved an asymptotic formula for it. In the present paper, first, we consider the Volterra integral equation for the remainder term in the asymptotic formula for the associated Euler totient function. Secondly, we solve the Volterra integral equation and we split the error term in the asymptotic formula for the associated Euler totient function into two summands called arithmetic and analytic part respectively.

\section{Introduction}
{\K} obtained a decomposition for the remainder term in the asymptotic formula for a generalization of the Euler totient function ({see~\cite{kac and wie}) : For a non-principal real Dirichlet character $\chi \hspace{0.1cm} (\text{mod }q), q>2,\text{ let } \varphi(n,\chi)$ denote the \textit{twisted Euler $\varphi$-function} 
\begin{equation*}
\varphi(n,\chi) = n\prod_{p|n} \left( 1 - \frac{\chi(p)}{p} \right).   \tag{1.1}
\end{equation*}
\footnote[0]
{\textit{2010 Mathematics Subject Classification. Primary}
 45D05;11A25;11N37

\textit{Key words and phrases}
: Volterra integral equation of second type,
the remainder term in the asymptotic formula,
twisted Euler $\varphi$-function,
associated Euler totient function.
}
Let
\vspace{-0.1cm}
\begin{equation*}
   E(x,\chi) = \sum_{n \leq x} \varphi(n,\chi) - \frac{x^2}{2L(2,\chi)}   \tag{1.2}
\end{equation*}and
\begin{equation*}
   E_1 (x,\chi) = \begin{cases}
                         E(x,\chi) \quad (x \notin \mathbb{N}),
                         \\
                         \frac{1}{2}(E(x-0,\chi)+E(x+0,\chi)) \quad (\text{otherwise})   \tag{1.3}
                      \end{cases}
\end{equation*}be the corresponding error terms. Here, as usual, $L(s,\chi)$ denotes the Dirichlet $L$-function associated to $\chi$. It is easy to see that $E(x,\chi) = O(x\log x)$ for $x \geq 2$. Hence $x^2 / (2L(2,\chi))$ is the main term in (1.2). Let $s(x)$ be the saw-tooth function
\vspace{-0.1cm}
\begin{equation*}
   s(x) = \begin{cases}
               0 \quad (x \in \mathbb{Z}),            
               \\
               \frac{1}{2} - \{ x \} \quad (\text{otherwise}),   \tag{1.4}
            \end{cases}
\end{equation*}where $\{ x \} = x - [x]$ is the fractional part of a real number $x$. We write for $x \geq 0$
\begin{gather*}
   f(x,\chi) = \sum_{n=1}^\infty \frac{\mu(n)\chi(n)}{n}s\left( \frac{x}{n} \right),   \quad \tag{1.5}
   \\
   g(x,\chi) = \sum_{n=1}^\infty \mu(n)\chi(n) \left\{ \frac{x}{n} \right\}\left( \left\{ \frac{x}{n} \right\} - 1 \right),  \quad  \tag{1.6}
\end{gather*}where $\mu(n)$ denotes the M\"{o}bious function. {\K} considered the Volterra integral equation of second type for (1.3) and solved it. 
\begin{theorem}[Theorem 1.1 in ~\cite{kac and wie}]
The solution of the following Volterra integral equation of second type
\begin{equation*}
   F(x,\chi) - \int_{0}^\infty K(x,t)F(t,\chi)dt = E_1 (x,\chi) \quad (x \geq 0),   \tag{1.7}
\end{equation*}where 
\begin{equation*}
   K(x,t) = \begin{cases}
                 1/t \quad (0 < t \leq x),
                 \\
                  0  \quad (0 \leq x < t), 
              \end{cases}
\end{equation*}is the function
\begin{equation*}
   F(x,\chi) = (f(x,\chi)+A)x,   \tag{1.8}
\end{equation*}where $A$ is an arbitrary constant.
\end{theorem}

(Probably, the term $Ax$ is missing to give the general solution as noted in ~\cite{I}.) Also, {\K} splitted (1.3) into two summands as follows :
\begin{theorem}[Theorem 1.2 in ~\cite{kac and wie}]
For $x \geq 0$ we have
\begin{equation*}
   E_1 (x,\chi) = xf(x,\chi) + \frac{1}{2}g(x,\chi).  \tag{1.9}
\end{equation*}By (1.9), $E_1 (x,\chi)$ can be splitted as follows :
\begin{equation*}
   E_1 (x,\chi) = E^{\text{AR}}(x,\chi) + E^{\text{AN}}(x,\chi) ,   \tag{1.10}
\end{equation*}where
\begin{equation*}
   E^{\text{AR}}(x,\chi) = xf(x,\chi) \quad \text{and} \quad E^{\text{AN}}(x,\chi) = \frac{1}{2}g(x,\chi)   \tag{1.11}
\end{equation*}with $f(x,\chi)$ and $g(x,\chi)$ given by (1.5) and (1.6) respectively.
We call $E^{\text{AR}}(x,\chi)$ and  $E^{\text{AN}}(x,\chi)$ the arithmetic part and analytic part of $E_1 (x,\chi)$ respectively.
\end{theorem}
 
{\ka} defined the associated Euler totient function for the generalized $L$-functions including the Riemann zeta function, Dirichlet $L$-functions and obtained the asymptotic formula (see ~\cite{Kac}) : By a polynomial Euler product we mean a function $F(s)$ of a complex variable $s=\sigma+it$ which for $\sigma>1$ is defined by a product of the form
\begin{equation*}
   F(s) = \prod_{p} F_p (s) = \prod_{p}\prod_{j=1}^d \left( 1-\frac{\alpha_j (p)}{ p^s } \right)^{-1},   \tag{1.12}
\end{equation*}where $p$ runs over primes and $|\alpha_j (p)| \leq 1$ for all $p$ and $1\leq j \leq d$. We assume that $d$ is chosen as small as possible, i.e. that there exists at least one prime number $p_0$ such that
\[ \displaystyle \prod_{j=1}^d \alpha_j ( p_0 ) \neq 0. \] Then $d$ is called the \textit{Euler degree} of $F$. For $F$ in (1.12) we define \textit{the associated Euler totient function} as follows :
\begin{equation}
  \varphi(n,F) = n\prod_{p|n} F_p (1)^{-1}  \quad (n \in \mathbb{N}).    \tag{1.13}
\end{equation}Let
\begin{gather*}
   \gamma(p) = p\left( 1-\frac{1}{F_p (1)} \right),   \tag{1.14}
   \\
   C(F) = \frac{1}{2}\prod_p \left( 1 - \frac{\gamma(p)}{p^2} \right).   \tag{1.15}
\end{gather*}{\ka} obtained the asymptotic formula for the error term in the asymptotic formula for (1.13).  
\begin{theorem}[Theorem 1.1 in ~\cite{Kac}]
For a polynomial Euler product $F$ of degree $d$ and $x \geq 1$ we have
\begin{equation*}
   \sum_{n \leq x} \varphi(n,F) = C(F)x^2 + O(x(\log 2x)^d).   \tag{1.16}
\end{equation*} 
\end{theorem}
Let us put
\begin{equation*}
   E(x,F) = \sum_{n \leq x} \varphi(n,F) - C(F)x^2   \tag{1.17}
\end{equation*}
and
\begin{equation}
   \alpha(n) = \mu(n)\prod_{p|n}\gamma(p),   \tag{1.18}
\end{equation}
where $\gamma(p)$ is defined by (1.14). The main aim of the present paper is to consider the Volterra integral equation of second type associated with $\varphi(n,F)$, and to prove the results similar to Theorem 1.1 and 1.2.

\section{Main theorems}
For a polynomial Euler product $F$ of degree $d$, let
\begin{equation*}
   E_2 (x,F)
   :=\begin{cases}
   E(x,F)  &(x \notin \mathbb{N})\\
   \frac{1}{2}(E(x-0,F)+E(x+0,F))  &(\textit{otherwise}).   \tag{2.1}
   \end{cases}\
\end{equation*} be the corresponding error terms. As in ~\cite{kac and wie}, we consider the following Volterra integral equation of second type for (2.1) as follows :
\begin{equation}
   F_1 (x,F)-\int_{0}^x F_1 (t,F)\frac{dt}{t} = E_2(x,F) \quad (x \geq 0),   \tag{2.2}
\end{equation}where $F_1 (x,F)$ is the unknown function. For every $x \geq 0$, let
\begin{equation} 
      f_1 (x,F) = \sum_{n=1}^\infty \frac{\alpha(n)}{n}s\left( \frac{x}{n} \right),  \tag{2.3}
\end{equation}where $s(x)$ is  the same as in (1.4). When $x$ is a positive integer, the following fact holds for $f_1 (x,F)$.
\begin{fact}
For positive integer $N$,
\begin{equation*}
   f_1 (N,F) = \frac{1}{2}(f_1 (N+0, F) + f_1 (N-0, F)).  \tag{2.4}
\end{equation*}
\end{fact}
\textit{Proof of Fact 2.1.} Let $N$ be a positive integer. By elementary calculations, we have
\begin{gather*}
   f_1 (N+0,F) = \frac{1}{2}\sum_{\substack{n\leq N+1\\ n|N}} \frac{\alpha(n)}{n}+\sum_{\substack{n\leq N+1 \\ n\nmid N}} \frac{\alpha(n)}{n}\left( \frac{1}{2}-\left\{ \frac{N+0}{n} \right\} \right)+\sum_{n>N+1}\frac{\alpha(n)}{n}\left( \frac{1}{2} - \frac{N}{n} \right), \tag{2.5}
\\
   f_1 (N-0,F) = -\frac{1}{2}\sum_{\substack{n\leq N+1\\ n|N}} \frac{\alpha(n)}{n}+\sum_{\substack{n\leq N+1 \\ n\nmid N}} \frac{\alpha(n)}{n}\left( \frac{1}{2}-\left\{ \frac{N-0}{n} \right\} \right)+\sum_{n>N+1}\frac{\alpha(n)}{n}\left( \frac{1}{2} - \frac{N}{n} \right).  \tag{2.6}
\end{gather*}Since
\[ \left\{ \frac{N+0}{n} \right\} + \left\{ \frac{N-0}{n} \right\} = 2\left\{ \frac{N+0}{n} \right\} \]for $n$ which does not divide a positive integer $N$, adding (2.5) and (2.6) we have
\begin{equation*}
   \frac{1}{2}(f_1 (N+0,F) + f_1 (N-0,F)) = \sum_{\substack{n\leq N+1 \\ n\nmid N}} \frac{\alpha(n)}{n}\left( \frac{1}{2}-\left\{ \frac{N+0}{n} \right\} \right)+\sum_{n>N+1}\frac{\alpha(n)}{n}\left( \frac{1}{2} - \frac{N}{n} \right).   \tag{2.7}
\end{equation*} By (2.3), the right-hand side of (2.7) corresponds to $f_1 (N,F)$.  $\square$ 

Moreover, to assure the convergence of the series (2.3), we assume that the series
\begin{equation*}
   \sum_{n=1}^\infty \frac{\alpha(n)}{n}  \tag{2.8}
\end{equation*}converges, where $\alpha(n)$ is the same in (1.18).
\begin{theorem}
For every complex number $A$, the function
\begin{equation}
   F_1 (x,F) = ( f_1 (x,F) + A)x  \quad (x \geq 0),    \tag{2.9}
\end{equation}is a solution of the integral equation (2.2) and these exhaust all solutions of (2.2).
\end{theorem}
As usual, in case we say a function $F_1 (x,F)$ is a solution of (2.2), we assume that the integral in (2.2) exists in the sense that the limit
\begin{equation*}
   \lim_{\epsilon \to 0+} \int_{\epsilon}^x |F_1 (t,F)|\frac{dt}{t}  \tag{2.10}  
\end{equation*}exists. We use the same convention throughout this paper. Also, the meaning of the integral in ~\cite{I} should also be interpreted in this sense. For every $x\geq0$, let
\begin{equation*}
   g_1 (x,F) = \sum_{n=1}^\infty \alpha(n)\left\{ \frac{x}{n} \right\} \left( \left\{ \frac{x}{n} \right\} -1 \right).  \tag{2.11}
\end{equation*}
\begin{theorem}
For $x \geq 1$ we have
\begin{equation*}
  E_2 (x,F) = x f_1 (x,F)+ \frac{1}{2}g_1 (x,F).   \tag{2.12}
\end{equation*}
\end{theorem}
We split $E_2 (x,F)$ into the \textit{arithmetic part} and the \textit{analytic part} as follows :
\begin{equation*}
  E_2 (x,F) = E^{\text{AR}}(x,F) + E^{\text{AN}}(x,F),   \tag{2.13}
\end{equation*}where
\begin{equation*}
   E^{\text{AR}}(x,F) = x f_1 (x,F) \quad \text{and} \quad E^{\text{AN}}(x,F) = \frac{1}{2}g_1 (x,F).   \tag{2.14}
\end{equation*}

\section{Remarks and auxiliary lemmas}
We prepare some remarks and auxiliary lemmas.
\begin{remark}[P33 in ~\cite{Kac}]
For every positive $\epsilon$, $\alpha(n) \ll n^\epsilon$. Hence the series
\begin{equation*}
   \sum_{n=1}^\infty \frac{\alpha(n)}{n^s}   \tag{3.1}
\end{equation*}absolutely converges for $\sigma>1$. Since $\alpha(n)$ is multiplicative by (1.18), we have 
\begin{equation*}
   \sum_{n=1}^\infty \frac{\alpha(n)}{n^2} = 2C(F).   \tag{3.2}
\end{equation*}
\end{remark}

\begin{remark}[Lemma 2.2 in ~\cite{Kac}]
The series
\begin{equation*}
   \sum_{n=1}^\infty \frac{\varphi (n,F)}{n^s}  \tag{3.3}
\end{equation*} converges absolutely for $\sigma > 2$ and in this half-plane we have
\begin{equation*}
    \sum_{n=1}^\infty \frac{\varphi(n,F)}{n^s} = \zeta(s-1)\sum_{n=1}^\infty \frac{\alpha(n)}{n^s}. \tag{3.4}
\end{equation*} In particular,
\begin{equation*}
   \varphi(n,F) = n\sum_{m|n} \frac{\alpha(m)}{m}.   \tag{3.5}
\end{equation*}
\end{remark}

We define the auxiliary function for $x \geq 0$ by
\begin{equation*}
   R(x,F) = E_2 (x,F) - xf_1  (x,F).   \tag{3.6}
\end{equation*}
\begin{lemma}
For all positive $x$,
 \begin{equation*}
   R(x, F) = -\int_{0}^x f_1 (t,F)dt.  \tag{3.7}
 \end{equation*}
\end{lemma}
\textit{Proof.} We can prove that $R(x,F)$ is a continuous function in the same way as in Lemma 1 of ~\cite{I}. For positive $x$ which is not an integer, take derivatives of the both sides of (3.6). Since $x$ is not a positive integer, we have $E_{2} ^\prime(x,F) = E^\prime (x,F) = -2C(F)x$. Therefore we have
\begin{equation*}
   R^\prime (x,F)
    = -2C(F)x - f_1 (x,F) - xf_{1}^\prime (x,F).
\end{equation*}Since $x$ is positive and not an integer, we have $\{x /n \}^\prime = 1/n$ (see ~\cite{Kac and Wie}, P2691). Considering the hypothesis on the series (2.8), Remark 3.1, and the fact that $x$ is positive and not an integer, differentiating term by term we obtain
          \begin{align*}
             \frac{d}{dx}\sum_{n=1}^\infty \frac{\alpha(n)}{n}s\left( \frac{x}{n} \right)
             &= \sum_{n=1}^\infty \frac{\alpha(n)}{n} \frac{d}{dx} \left( \frac{1}{2} - \left\{ \frac{x}{n} \right\} \right)
             \\
             &= -2C(F). 
          \end{align*}Consequently, we have 
\[ R^\prime (x,F) = -f_1 (x,F) \] for  $x$ which is positive and not an integer. Since $R(0,F) = 0$ by (2.3) and $R(x,F)$ is continuous for all positive $x$, we have (3.7) for all positive $x$. $\square$

\begin{lemma}Let $G$ be a complex-valued function defined on $[0, \infty)$ satisfying
                         \begin{equation}
                            \int_{0}^x |G(t)|\frac{dt}{t} < +\infty   \tag{3.8}
                         \end{equation} and the integral equation
                         \begin{equation}
                            G(x) - \int_{0}^x G(t)\frac{dt}{t} = 0   \tag{3.9}
                         \end{equation}for all $x \geq 0$. Then we have
                         \begin{equation}
                            G(x) = Ax    \tag{3.10}
                        \end{equation}for some complex number $A$.
\end{lemma}
\textit{Proof.} This is Lemma 2 in ~\cite{I}. $\square$

\section{Proof of main theorems.}
First we prove Theorem 2.1 for $x$ which is positive and not an integer. Let a function $F_1 (x,F)$ be a solution of the Volterra integral equation of second type (2.2) satisfying the condition (2.10). Using Lemma 3.3, from (3.6) we have 
\begin{equation*}
E_2 (x,F) - xf_1 (x,F) = -\int_{0}^x f_1 (t,F)dt.   \tag{4.1}
\end{equation*} Since $x$ is positive and not an integer, $E_2 (x,F) = E(x,F)$.  By (2.2), we have
         \begin{equation*}
             \int_{0}^x (F_1 (t,F)-tf_1 (t,F))\frac{dt}{t} = F_1 (x,F) - xf_1 (x,F).
         \end{equation*}Using Lemma 3.4, we have the solution \[ F_1 (x,F) = (f_1 (x,F) + A)x. \] Conversely, if we assume that $F_1 (x,F)$ is a function of type (2.9). Then,  by (3.6) and (3.7),
\begin{align*}
   F_1 (x,F) - \int_{0}^x F_1 (t,F) \frac{dt}{t}
   &= xf_1 (x,F) - \int_{0}^x f_1 (t,F)dt     
    \\
    &= xf_1 (x,F) + R(x,F)
   \\
   &= E (x,F).
\end{align*} Therefore, the function $F_1 (x,F)$ of type (2.5) is a solution of (2.2) for positive and not an integer $x$. Also, the function $f_1 (x,F)$ is a locally bounded. In fact, by the hypothesis (2.8) and (3.2)
\begin{align*}
   f_1 (x,F)
   &=  \sum_{n=1}^\infty \frac{\alpha(n)}{n}s\left( \frac{x}{n} \right)
   \\
   &= \sum_{n=1}^\infty \frac{\alpha(n)}{n} \left( \frac{1}{2} - \left\{ \frac{x}{n} \right\} \right)
   \\
   &= \frac{1}{2}\sum_{n=1}^\infty \frac{\alpha(n)}{n} -  \sum_{n=1}^\infty \frac{\alpha(n)}{n}\left\{ \frac{x}{n} \right\}
   \\
   &= \frac{1}{2}\sum_{n=1}^\infty \frac{\alpha(n)}{n} - \sum_{n=1}^\infty \frac{\alpha(n)}{n} \left( \frac{x}{n} - \left[ \frac{x}{n} \right] \right)
   \\
   &= \frac{1}{2}\sum_{n=1}^\infty \frac{\alpha(n)}{n} - 2C(F)x + \sum_{n \leq x} \frac{\alpha(n)}{n}\left[ \frac{x}{n} \right].         
\end{align*}It is clear that the function $F_1 (x,F)$ satisfies the condition (2.10).

Next we prove Theorem 2.1 for $x$ which is a positive integer. Let a function $F_1 (x,F)$ be the solution of the equation (2.2) satisfying the condition (2.10). Since $x$ is a positive integer, $E_2 (x,F) = \frac{1}{2}(E(x+0,F) + E(x-0,F))$. By continuity of $R(x,F)$ for all positive $x$ and (2.4), we obtain
\begin{align*}
   F_1 (x,F) - \int_{0}^x F_1 (t,F)\frac{dt}{t}
   &= \frac{1}{2}(E(x+0,F) + E(x-0,F))
   \\
   &= R(x,F) + xf_1 (x,F).  \tag{4.2} 
\end{align*}Using (3.7) and Lemma 3.4, we see that the function (2.9) is the solution of (2.2). Conversely, if we assume that $F_1 (x,F)$ is a function of type (2.9). By substituting it into the left hand side of (2.2), we have
\begin{equation*}
            F_1 (x,F)- \int_{0}^x F_1 (t,F)\frac{dt}{t}  
            = f_1 (x,F) - \int_{0}^x f_1 (t,F)dt.
\end{equation*}Using (2.2), 
\begin{gather*}
            E(x+0, F)
            = F_1 (x+0, F) - \int_{0}^x F_1 (t,F)\frac{dt}{t},  \tag{4.3}
            \\
            E(x-0, F)
            = F_1 (x-0, F) - \int_{0}^x F_1 (t,F)\frac{dt}{t}.  \tag{4.4}
\end{gather*}By (4.2), (4.3), (2.9) and (2.4), we have
\begin{equation*}
\frac{1}{2}(E(x+0,F) + E(x-0,F)) = f_1 (x,F) - \int_{0}^x f_1 (t,F)dt. 
 \end{equation*}Since $x$ is a positive integer, the left hand side corresponds to $E_2 (x,F)$. Therefore, Theorem 2.1 also holds for all positive integer $x$.   $\square$

Let us prove Theorem 2.2. By lemma 3.3 it is enough to show that for $x \geq 1$ we have
    \begin{equation*}
       \int_{0}^x f_1 (t,F)dt = -\frac{1}{2}g_1 (x,F).   \tag{4.5}
    \end{equation*} This can be done as follows : Recalling Lemma 3.3, we have
\begin{equation*}
   R(x,F) =  -\frac{x}{2}\left( \sum_{n=1}^\infty \frac{\alpha(n)}{n} \right) + \frac{1}{2}\sum_{n=1}^\infty \alpha(n)\frac{x^2}{n^2}  - \sum_{n \leq x} \frac{\alpha(n)}{n}\int_{n}^x \left( \frac{t}{n} - \left\{ \frac{t}{n} \right\} \right)dt.   \tag{4.6}
\end{equation*} For  $x>0$, 
\begin{equation*}
   \int_0^x \{ t \}dt = \frac{1}{2}\{ x \}^2 + \frac{1}{2}[x]
\end{equation*}(see ~\cite{Kac and Wie}, P2692) and hence we have
\begin{equation*}
   \int_{n}^x \left\{ \frac{t}{n} \right\}dt = \frac{n}{2}\left( \left\{ \frac{x}{n} \right\}^2 + \left[ \frac{x}{n} \right] - 1 \right).   \tag{4.7}  
\end{equation*} By substituting it into (4.6), we have
\begin{align*}
   R(x,F) 
   &= -\frac{1}{2}\sum_{n=1}^\infty \alpha(n)\frac{x}{n} + \frac{1}{2}\sum_{n=1}^\infty \alpha(n)\frac{x^2}{n^2} - \frac{1}{2}\sum_{n \leq x} \alpha(n) \left( \frac{x}{n} - \left\{ \frac{x}{n} \right\} - 1 \right)\left[ \frac{x}{n} \right]
   \\
   &=  -\frac{1}{2}\sum_{n=1}^\infty \alpha(n)\frac{x}{n} + \frac{1}{2}\sum_{n=1}^\infty \alpha(n)\frac{x^2}{n^2} - \frac{1}{2}\sum_{n=1}^\infty \alpha(n) \left( \frac{x}{n} - \left\{ \frac{x}{n} \right\} - 1 \right)\left[ \frac{x}{n} \right]
   \\
   &= -\frac{1}{2}\sum_{n=1}^\infty \alpha(n)\left\{ \frac{x}{n} \right\} + \frac{1}{2} \sum_{n=1}^\infty \alpha(n) \left\{ \frac{x}{n} \right\}^2
   \\
   &= \frac{1}{2}g_1 (x,F).    \tag{4.8}
\end{align*} The proof is complete. $\square$

\vspace{0.25cm}

Theorems 2.1 and 2.2 are generalizations of Theorems 1.1 and 1.2. This can be seen as follows : If $F$ is the Dirichlet $L$-function $L(s,\chi)$ in (1.13), then by (1.14) the associated Euler totient function $\varphi(n,F)$ corresponds to the twisted Euler $\varphi$-function $\varphi(n,\chi)$. Since the Euler degree of $L(s,\chi)$ equals to $1$, we have $\gamma(p) = \chi(p)$ in (1.14). By (1.18), we have $\alpha(n) = \mu(n)\chi(n)$. Therefore we have
\begin{gather*}
   f_1 (x,F) = \sum_{n=1}^\infty \frac{\mu(n)\chi(n)}{n}s\left( \frac{x}{n} \right) = f(x,\chi),   \tag{4.9}
   \\
   g_1 (x,F) = \sum_{n=1}^\infty \mu(n)\chi(n)\left\{ \frac{x}{n} \right\} \left( \left\{ \frac{x}{n} \right\} -1 \right) = g(x,\chi).   \tag{4.10}
\end{gather*} Hence, (4.9) and (4.10) correspond to (1.5) and (1.6) respectively.

\vspace{0.25cm}

\textbf{Acknowledgments.}   The author expresses his sincere gratitude to Prof. Kohji Matsumoto.

\bigskip
　　　　　　　　　　　　　　Hideto Iwata 
                                                  
　　　　　　　　　　　　　　Graduate School of Mathematics

　　　　　　　　　　　　　　 Nagoya University
                                                   
　　　　　　　　　　　　　　 Furocho, Chikusa-ku,
                                                   
　　　　　　　　　　　　　　 Nagoya, 464-8602, Japan.

　　　　　　　　　　　　　　 \small{e-mail:d18001q@math.nagoya-u.ac.jp}

\begin{thebibliography}{99}
\bibitem{I}
H. Iwata, 
\newblock On the solution of the Volterra integral equation of second  type for the error team in an asymptotic formula for arithmetic functions, 
\newblock  Advanced Studies : Euro-Tbilisi Mathematical Journal (to appear).

\bibitem{Kac and Wie} 
J. Kaczorowski, K. Wiertelak, 
\newblock Oscillations of the remainder term related to the Euler totient function, 
\newblock J. Number Theory \textbf{130} (2010) 2683-2700. 

\bibitem{kac and wie}
J. Kaczorowski, K. Wiertelak, 
On the sum of the twisted Euler function,  
\newblock Int. J. Number Theory \textbf{8} (7) (2012) 1741-1761. 

\bibitem{Kac}
J.Kaczorowski,
On a generalization on the Euler totient function,
\newblock Monatsh Math (2013) \textbf{170} : 27-48.

\end{thebibliography}
\end{document}